\documentclass[12pt]{amsart}
\usepackage{amssymb}
\input psfig

\newcommand{\z}{{\mathbb Z}}

\newcommand{\zp}{\z_p}
\newcommand{\zpn}{\z/p^n\z}
\newcommand{\zpnpone}{\z/p^{n+1}\z}
\newcommand{\zpnmone}{\z/p^{n-1}\z}
\newcommand{\zpp}{\z/p\z}
\newcommand{\qp}{{\mathbb Q}_p}
\newcommand{\ord}{{\text{ord}}_p}
\newcommand{\tth}{^{\text{th}}}
\newcommand{\xt}{\tilde{x}}
\newcommand{\st}{\tilde{\sigma}}
\DeclareMathOperator{\OO}{O}
\newcommand{\PP}{{\mathbb P}}
\newcommand{\opshrink}[2]{\mskip -#1mu minus -#1mu #2 \mskip -#1mu minus -#1mu}
\newcommand{\mycirc}{\opshrink{2}{\ncirc}}
\newcommand{\ncirc}{
   \mathbin{\vcenter{\hbox{\setlength{\unitlength}{.6ex}%
               \begin{picture}(2.0,0)(-1.0,0)
               \put(0,0){\circle{1.0}}
               \end{picture}}}}}
\newtheorem{lemma}{Lemma}
\newtheorem{corollary}{Corollary}
\newtheorem{proposition}{Proposition}

\begin{document}
\title[Polynomial Mappings mod $p^n$]
{On the Structure of Polynomial Mappings Modulo an Odd Prime Power}

\thanks{This paper is unchanged from the version circulated in August 1994.  See the
second author's Ph.D. thesis (Berkeley, 1996) for further results and for comments
relating this paper to the mathematical literature.}

\author{David L. desJardins}
\address{Google, Inc.,
         2400 Bayshore Parkway,
         Moutain View, CA 94043}
\email{david@desjardins.org}
\author{Michael E. Zieve}
\address{Center for Communications Research,
         29 Thanet Road,
         Princeton, NJ 08540}
\email{zieve@idaccr.org}


\begin{abstract}
Let $f(x) \in \z[x]$ be a polynomial with integer coefficients, let $n$
be a positive integer, and $p$ an odd prime.  Then the mapping $x
\mapsto f(x)$ sends $\zpn$ into $\zpn$.  We study the topological
structure of this mapping.
\end{abstract}

\date{August 1994}             

\maketitle


\section {Introduction}

Let $f(x) \in \z[x]$ be a polynomial with integer coefficients, let $n$
be a positive integer, and let $p$ be an odd prime.  Then the mapping $x
\mapsto f(x)$ sends $\zpn$ into
$\zpn$.  We shall study the structure of this mapping.  Since the
mapping $f\pmod{p^n}$ must project to a well-defined mapping
$f\pmod{p^{n-1}}$, only a certain class of mappings on $\zpn$ can arise
from polynomials.  But there are many more restrictions on which
mappings can occur than just the above observation---in
Section~\ref{lifting} we show that there is a certain linearity causing
one such restriction.  In later sections we take advantage of this
linearity to derive numerous results about the cycles of $f\pmod{p^n}$.
Our results give an algorithm which, for almost any given polynomial
$f$, finds the lengths of the cycles of $f\pmod{p^n}$ for all $n$,
usually very quickly.  Our results also indicate how to construct a
polynomial with any (possible) desired cycle structure mod~$p^n$.
Our methods also apply in much more general situations\footnote{For
instance, we can allow our polynomials to have coefficients in the
$p$-adic integers $\zp$, and in fact every argument we make will be
unchanged if we replace every symbol $\z$ by the symbol $\zp$.}; we will
briefly discuss this in Section~\ref{etc.}.


\section{Notation}

Henceforth, $p$ will denote a fixed odd prime, $f(x) \in \z [x]$ a fixed
polynomial, and $n$ a positive integer.  We denote by $f_n$ the mapping
$\zpn \rightarrow \zpn$ which sends $x \mapsto f(x) \pmod{p^n}$.  We let
$\sigma = (x_1,\dots, x_k)$ be a cycle of $f_n$ of length $k$; that is,
$f_n(x_1)=x_2$, $f_n(x_i)=x_{i+1}$, and $f_n(x_k)=x_1$.  (We view the
$x_i$ as integers lying in the appropriate classes (mod~$p^n$).)
Finally, $g=f^k$ is the $k\tth$ iterate of $f$.


\section{Cycle lifting} 
\label{lifting}

In this section we examine the structure of $f_{n+1}$ on the set of
points of $\zpnpone$ which are congruent mod $p^n$ to elements of
$\sigma$.  Let $X_i$ be the preimage of $x_i$ under the projection
$\zpnpone \rightarrow \zpn$; thus, $|X_i|=p$, and by the definition of
$\sigma$, $f_{n+1}(X_i) \subseteq X_{i+1}$.  For $g=f^k$, the $k\tth$
iterate of $f$, we have $g_{n+1}(X_1) \subseteq X_1$.

\begin{figure}[h]
\centerline{\psfig{figure=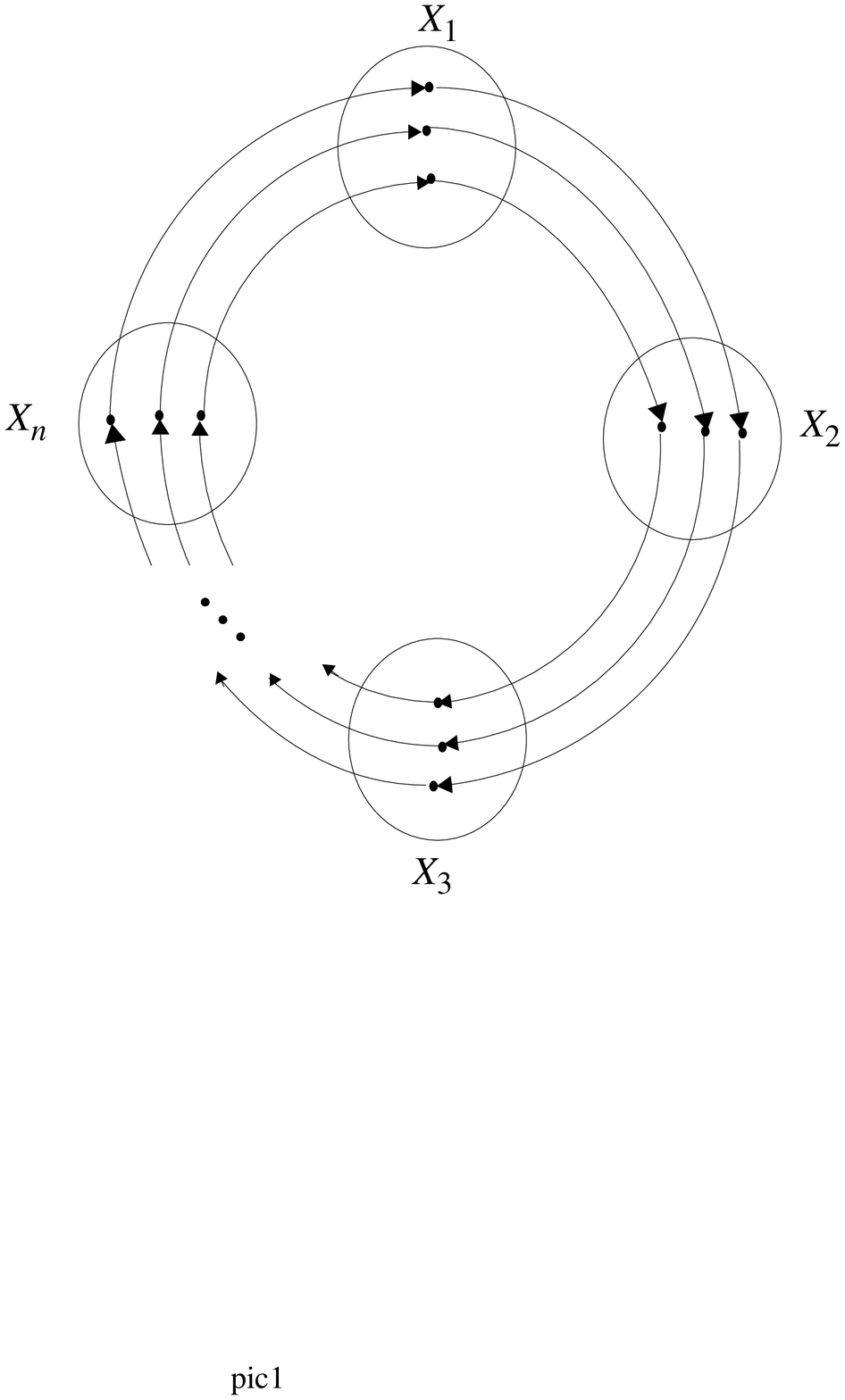,height=4in,width=4in,bblly=4in,bbury=10in,bbllx=1in,bburx=7in,clip=}}
\end{figure}

Now, let $X=X_1 \bigcup X_2 \bigcup \dots \bigcup X_k$; then $f_{n+1}(X)
\subseteq X$, and any cycle of $f_{n+1}$ in $X$ must have length
divisible by $k$.  We call these cycles the {\em lifts} of $\sigma$.
When we divide the lengths of these lifts by $k$, we get the cycle
lengths of $g_{n+1}$ in $X_1$.

We can define a bijection between $X_1$ and $\zpp$ by the rule
\begin{equation*}
x_1 + p^n t \longleftrightarrow t.
\end{equation*}
By Taylor's theorem for polynomials,\footnote{This says that $g(x+y)=\sum_{i=0}
^{\infty} y^i \frac{g^{(i)}(x)}{i!}$; note that the sum is finite, since all
terms with $i>\text{degree}\ (g)$ vanish, and also note that $\frac{g^{(i)}(x)}
{i!}$ is an integer.}
\begin{align*}
g(x_1 + p^n t) & \equiv g(x_1) + p^n t g'(x_1) \pmod{p^{2n}} \\
&\equiv x_1 + p^n \left (\frac{g(x_1) -x_1}{p^n} \right ) + p^n g'(x_1) t
\pmod{p^{2n}} \\
&\equiv x_1 + p^n b_n + p^n a_n t \pmod{p^{2n}}
\end{align*}
where we define $a_n = g'(x_1)$ and $b_n = (g(x_1) - x_1)/p^n$.  (Note that
$a_n$ and $b_n$ are defined over $\z$.)  Thus, if we define the map
$\Phi:\zpp \rightarrow \zpp$ to be induced by restricting $g_{n+1}$ to
$X_1$ and applying the above bijection, then $\Phi(t) = b_n + a_n t$.
The linearity of this map is the key to what follows.

Note that:
\begin{enumerate}
\item 
If $a_n \equiv 1 \pmod{p}$ and $b_n \not\equiv 0 \pmod{p}$, then $\Phi$
consists of a single cycle of length $p$, so that $f_{n+1}$ restricted
to $X$ consists of a single cycle of length $pk$.  In this case we say
that $\sigma$ {\em grows}.
\item
If $a_n \equiv 1 \pmod{p}$ and $b_n \equiv 0 \pmod{p}$, then $\Phi$ is
the identity, so $f_{n+1}$ restricted to $X$ consists of $p$ cycles,
each of length $k$.  In this case we say that $\sigma$ {\em splits}.
\item
If $a_n \equiv 0 \pmod{p}$, then $\Phi$ is constant, so $f_{n+1}$ on $X$
contains one $k$-cycle, and the remaining points of $X$ are mapped into
this cycle by $f^k$.  In this case we say that $\sigma$ {\em grows
tails}.
\item
If $a_n \not\equiv 0,1 \pmod{p}$, then $\Phi$ is a permutation, and
$\Phi^\ell$, the $\ell\tth$ iterate of $\Phi$, sends 
\begin{align*}
t&\rightarrow (b_n+a_nb_n+a_n^2b_n+\dots+a_n^{\ell-1}b_n)
 + a_n^\ell t\\ 
&= b_n(a_n^\ell-1)/(a_n -1) +a_n^\ell t,
\end{align*}
so $\Phi^\ell(t)-t=(t+b_n/(a_n -1))(a_n^\ell -1)$.  Thus, $\Phi$ has a
single fixed point, namely $t=-b_n/(a_n-1)$, and the remaining points
of $X_1$ lie on cycles of length $d$, where $d$ is the order of $a_n$
in $(\zpp)^*$.  Thus, $f_{n+1}$ restricted to $X$ contains one
$k$-cycle and $(p-1)/d$ cycles of length $kd$.  In this case we say
that $\sigma$ {\em partially splits}.
\end{enumerate}

\begin{figure}[h]
\centerline{\psfig{figure=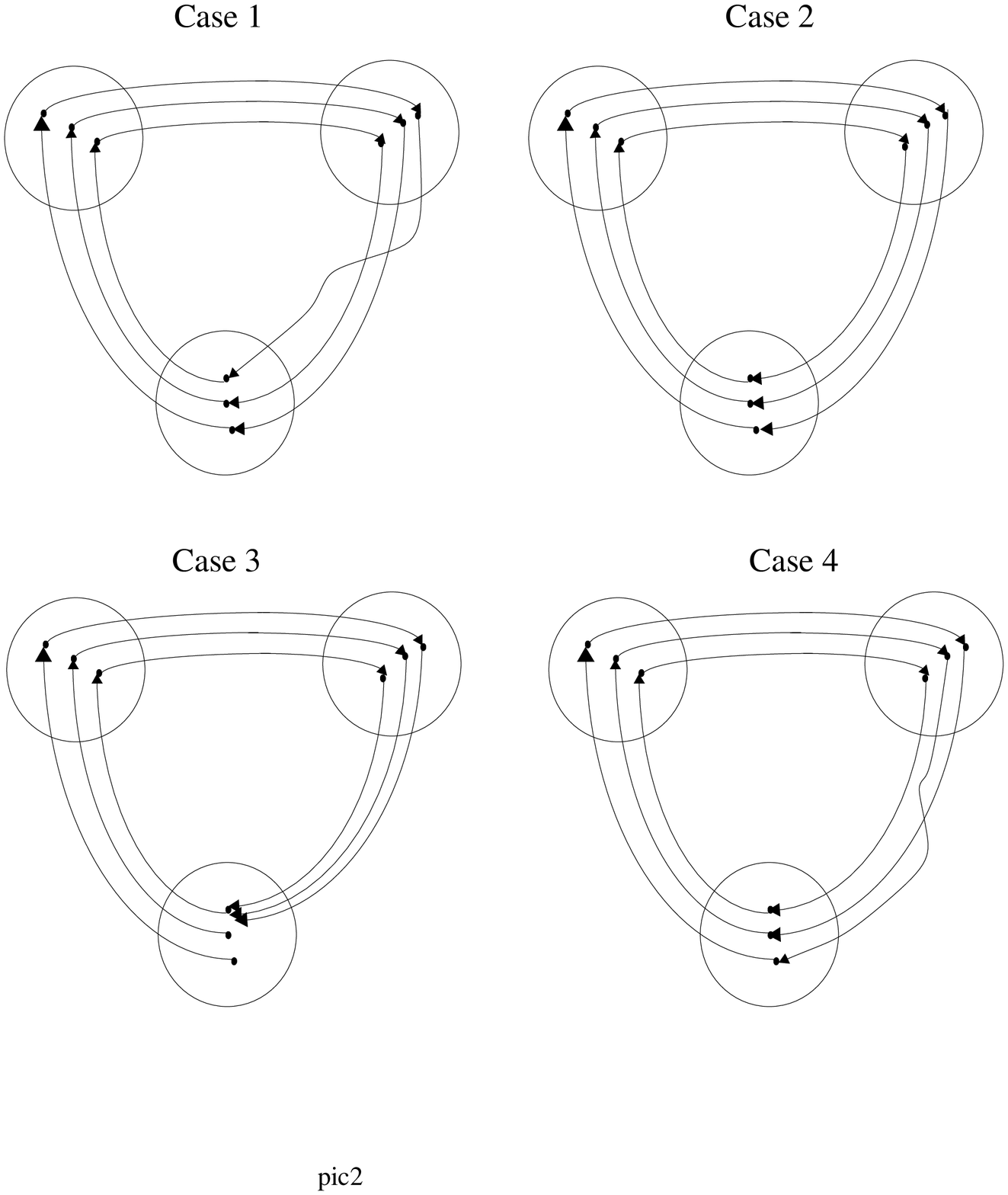,height=6in,width=5in,bblly=2in,bbury=11in,bbllx=0in,bburx=8in,clip=}}
\end{figure}

The above results already rule out many potential ways for $\sigma$ to
lift.  For instance, the lifts of $\sigma$ can only have two distinct
lengths.  If two lengths do occur, then one equals the length of
$\sigma$ and occurs only once; if only one length occurs, it is either
the length of $\sigma$ or $p$ times this length.

Before proceeding any further, we comment on the definitions of $a_n$
and $b_n$.  Our definition of $a_n$ and $b_n$ depends on three things:
the cycle $\sigma$, the choice of $x_1$ from the elements of $\sigma$,
and the integer $x_1$ chosen to represent the congruence class $x_1
\pmod{p^n}$.  However, to some extent $a_n$ and $b_n$ are independent of
these last two choices.  First,
\begin{equation*}
a_n= (f^k)'(x_1) = \prod_{i=0}^{k-1}f'(f^i(x_1))
\equiv \prod_{i=0}^{k-1}f'( x_i) \pmod{p^n},
\end{equation*}
so the class of $a_n \pmod{p^n}$ does not depend on the choices.
Secondly,
\begin{equation*}
g(x_1 +p^nz) -(x_1+p^nz) \equiv p^nb_n +p^nz(a_n -1)
\pmod{p^{2n}},
\end{equation*}
so replacing $x_1$ by $x_1+p^nz$ has the effect of
replacing $b_n$ with $b_n+z(a_n-1) \pmod{p^n}$.  Thus, for
$A=\min\{\ord(a_n-1), n\}$, the choice of the integer $x_1$ from the
congruence class $x_1\pmod{p^n}$ does not affect $b_n \pmod{p^A}$.
Finally,
\begin{equation*}
g(f(x_1))-f(x_1)=f(x_1+p^nb_n)-f(x_1) \equiv p^nb_nf'(x_1)
\pmod{p^{2n}},
\end{equation*}
and since $p\nmid f'(x_1)$, $\min\{\ord(b_n),n\}$ is independent of
the choice of which a particular element of $\sigma$ is called $x_1$.


\section{Relationships between $a$'s and $b$'s}
\label{relationships}

Let $\st=(\xt_1,\dots,\xt_{rk})$ be a lift of $\sigma$ to an $rk$-cycle
of $f_{n+1}$.  We will show that the manner in which $\sigma$ lifted
restricts how $\st$ can lift.  We may assume that $\xt_1 \equiv x_1
\pmod{p^n}$, and (as before, viewing $\xt_1$ as an integer) we write
$\xt_1=x_1 + p^nt$.  Then
\begin{align*}
a_{n+1} = (g^r)'(x_1 + p^nt) &\equiv (g^r)'(x_1) = \prod_{i=0}^{r-1} 
g'(g^i(x_1)) \\
&\equiv  g'(x_1)^r = a_n^r \pmod{p^n}.
\end{align*}
Now we apply this calculation:
\begin{enumerate}

\item
If $\sigma$ splits or grows, then $a_{n+1} \equiv 1^r \equiv 1
\pmod{p}$, so $\st$ either splits or grows.

\item 
If $\sigma$ partially splits, then its $k$-cycle lift also partially
splits (with the same $d$, since $a_{n+1}\equiv a_n \pmod{p}$ and so the
order of $a_{n+1}$ in $(\zpp)^*$ is the same as the order of $a_n$ in
$(\zpp)^*$), and its $kd$-cycle lifts either split or grow (since
$a_{n+1} \equiv a_n^d \equiv 1 \pmod{p}$).

\item
If $\sigma$ grows tails, then the single $k$-cycle lift $\st$ also has
$a_{n+1} \equiv a_n \equiv 0 \pmod{p}$, so it grows tails as well.
\end{enumerate}

We will need another basic calculation.  As before,
\begin{equation*}
g^r(x_1 + p^n t)
\equiv x_1 + p^n (ta_n^r + b_n (1+a_n +\dots +a_n^{r-1})) \pmod{p^{2n}},
\end{equation*}
so
\begin{align*}
p^{n+1}b_{n+1} &= g^r (x_1 + p^n t) - (x_1 + p^n t) \\
&\equiv  p^n \bigl(t(a_n^r -1) + b_n (1+a_n +\dots +a_n^{r-1})\bigr)
\pmod{p^{2n}}
\end{align*}
and therefore
\begin{equation*}
pb_{n+1} \equiv t(a_n^r -1) + b_n (1+a_n +\dots + a_n^{r-1}) \pmod{p^n}.
\end{equation*}


\section{Outline of goals} 
\label{philosophy}

Now that we have established the basic setup, we briefly pause to
discuss the general questions we are studying.  We have seen that the
cycle structure of $f_n$ greatly depends on that of $f_{n-1}$.  Thus, it
will be possible to obtain results which apply to the structure of $f_n$
for all $n$.  More precisely, we study an infinite tree which contains a
node for each cycle of $f_n$, for every $n\ge 0$, and where each node is
labeled with the length of the corresponding cycle.  The tree is defined
as follows: at the top level, level $0$, is a single node labeled with
1, the length of the single cycle of $f_0$.  At each lower level, level
$n$, there is a node for each cycle of $f_n$, labeled with its length,
and a node at level $n+1$ is a child of a node at level $n$ if it is a
lift of the corresponding cycle.

Here is an example of such a tree, for a polynomial with $p=3$:

\begin{figure}[h]
\centerline{\psfig{figure=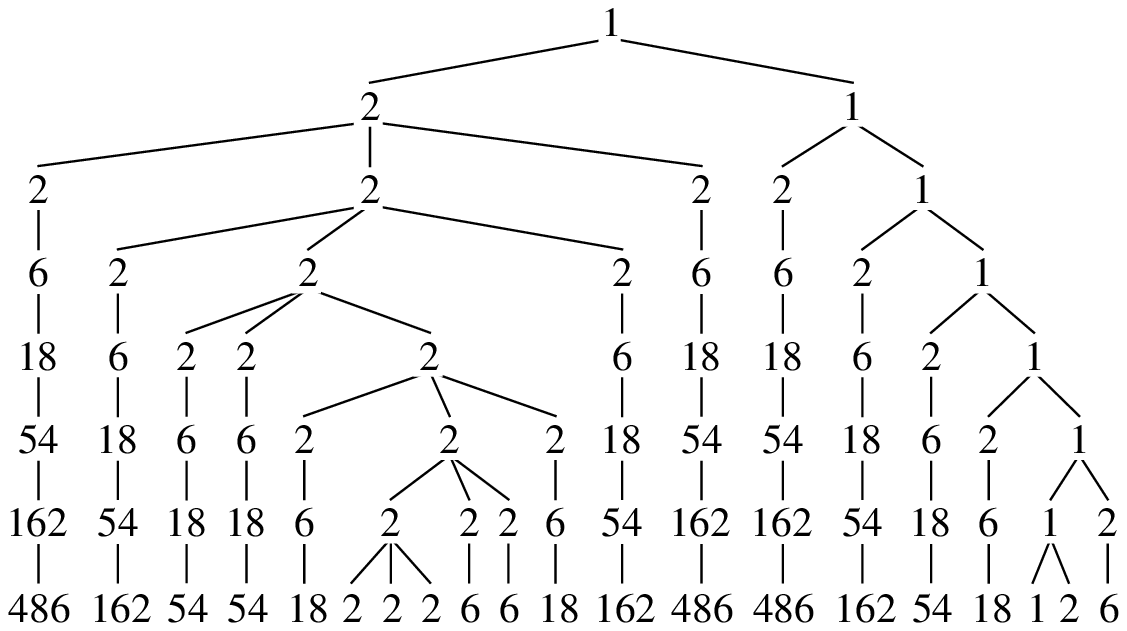,bblly=1in,bbury=4in,bbllx=1in,bburx=6in,clip=}}
\end{figure}

We would like to do the following:
\begin{enumerate}
\item Describe all trees that can occur.
\item Give a method for constructing a polynomial having a prescribed tree.
\item Give a method for determining the tree of a given polynomial.
\end{enumerate}
We will derive a number of results of the form: if a certain (finite)
part of the tree has a certain form, then this constrains the behavior
of another (possibly infinite) part of the tree.  For example, we will
show that whenever a cycle for some $f_n$ (with $n\geq 2$) grows, then
its lift grows, and the lift of that lift grows, and so on.

Results of this form severely restrict the class of trees which can
occur.  They are also useful for determining the tree of a given
polynomial; in fact, except for a certain pathological class of
functions, we will see that the first $n$ levels of the tree, for some
$n$, will determine the entire tree.  We have looked at thousands of
random polynomials of small degree, and in every case the first nine levels
were sufficient; usually five were enough, and it seemed that fewer levels
were needed for larger $p$.  The pathological cases appear to be quite
rare, since none arose randomly.  However, we do not believe that, in
the non-pathological cases, there is a bound on the number of levels of
the tree needed to determine the entire tree; large numbers of levels
should sometimes be necessary, but only very rarely.

Finally, these results help us construct polynomials having prescribed
trees.  As long as the tree is determined by its first $n$ levels, we
need only find a polynomial whose tree has those first $n$ levels; i.e.,
a polynomial having a certain structure mod $p^n$.

We will generally not study cycles which grow tails, except in
Section~\ref{tails}.  This case is easy to identify and distinguish,
because cycles which grow tails will only occur in subtrees rooted at
cycles mod $p$ which grow tails.


\section{Cycle structures}


\subsection{If $\sigma$ grows}

\label{grow}
Suppose $\sigma$ grows.  We showed above that $\st$ either splits or
grows.  From $a_n \equiv 1 \pmod{p}$, it follows that $a_n^p \equiv 1
\pmod{p^2}$, because
\begin{equation*}
\frac{a_n^p-1}{a_n-1} = a_n^{p-1}+\dots+a_n+1 \equiv 1+\dots+1+1
\equiv 0 \pmod{p}.
\end{equation*}
So, for $n\geq 2$,
\begin{equation*}
pb_{n+1} \equiv b_n (1+a_n +\dots +a_n^{p-1}) \pmod{p^2}.
\end{equation*}
If $a_n =1$, then $1+a_n +\dots +a_n^{p-1} =p$.  Otherwise, let $a_n =
1+ p^\gamma \delta$, where $p\nmid\delta$ and $\gamma \geq 1$.  Then
\begin{align*}
1+a_n +\dots +a_n^{p-1} = \frac{a_n^p -1}{a_n -1} &=
\frac{{p\choose 1}p^\gamma\delta + {p\choose 2}p^{2\gamma}\delta^2 +\cdots}
{p^\gamma\delta} \\
&={p\choose 1}+{p\choose 2}p^\gamma\delta + \cdots \\
&\equiv p \pmod{p^2}.  
\end{align*}
Thus, in either case $pb_{n+1} \equiv pb_n \pmod{p^2}$.  So, $b_{n+1}
\not\equiv 0 \pmod{p}$, so $\st$ grows.

Here we have shown that, for $n\geq 2$, whenever $\sigma$ grows, its
lift $\st$ also grows; it follows that the lift of $\st$ also grows, and
so on.  In Appendix~\ref{growth-appendix} we will show that, for $p>3$,
this result holds for $n=1$ as well.

In this case, the subtree rooted at $\sigma$ has the following
structure:

\begin{figure}[h]
\centerline{\psfig{figure=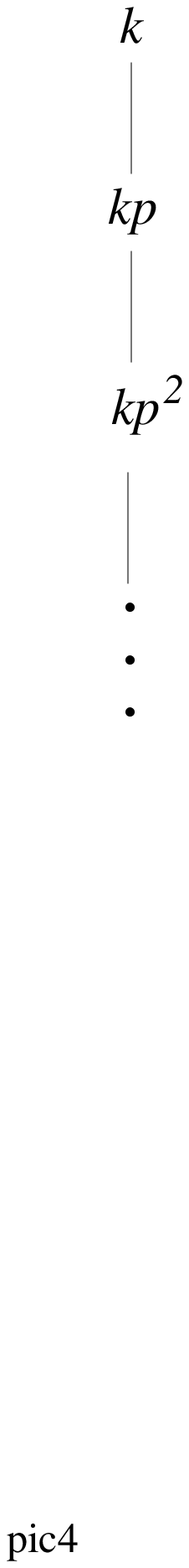,height=2.8in,bblly=5in,bbury=8.8in,bbllx=3in,bburx=4in,clip=}}
\end{figure}


\subsection{If $\sigma$ partially splits}

\label{partial}
Suppose that $\st$ is a $kd$-cycle lift of $\sigma$, where $1<d<p$;
thus, $d$ is the order of $a_n$ in $(\zpp)^*$.  Let $y$ be any element
of $\st$, and let $h=f^{kd}$.  In Corollary~\ref{useful} of
Section~\ref{improvements} we show that 
\begin{equation*}
\min\{\ord(h(y)-y)-n,nd\} =
\min\{\ord(a_{n+1}-1),nd\}.
\end{equation*}
In this section we note some implications of this result:
\begin{enumerate}
\item
If $e=\ord (a_{n+1} -1) <nd$, then $h(y)\equiv y\pmod p^{n+m}$ for $m\le
e$, but $h(y)\not\equiv y\pmod p^{n+e+1}$, so $\st$ splits $(e-1)$
times, and then the descendants of $\st$ at level $n+e$ grow.
\item
If $e=\ord(a_{n+1}-1)\geq nd$, then $\st$ splits $(nd-1)$ times, but we
do not know what happens to its descendants at level $n+nd$.
\end{enumerate}
Note that, if $e<n$, then every $kd$-cycle lift of $\sigma$ has the same
$e$ (since $a_{n+1} \equiv a_n^d \pmod{p^n}$), so they all behave the
same way.

In case 1, the subtree rooted at $\sigma$ has the following structure:
\pagebreak
\begin{figure}[h]
\centerline{\psfig{figure=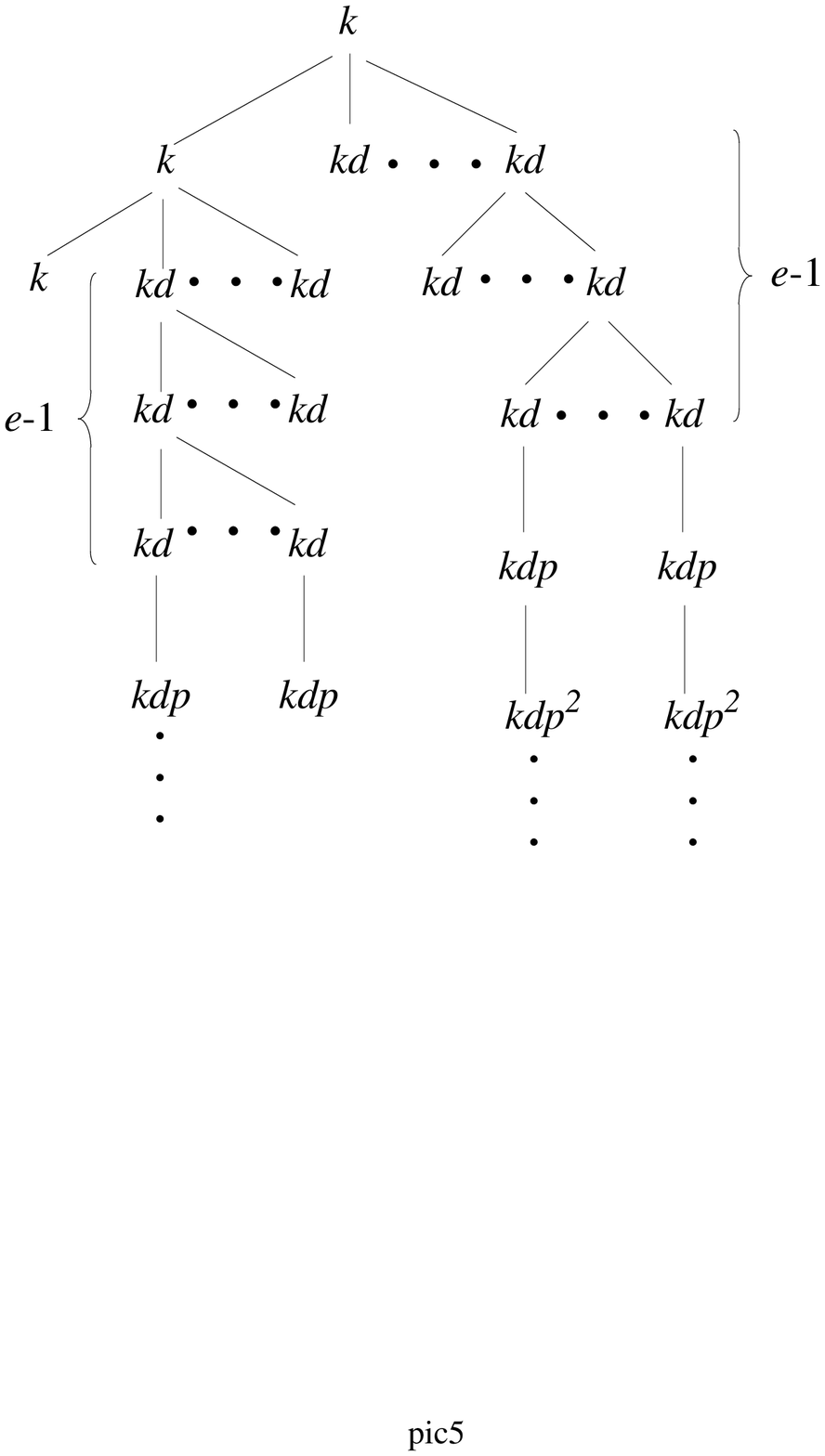,height=4.5in,width=5in,bblly=5in,bbury=11in,bbllx=0in,bburx=7in,clip=}}
\end{figure}


\subsection{If $\sigma$ splits}
\label{split}

Suppose that $\sigma$ splits.  Let $x$ be an element of $\sigma$, and
let $y=x+p^nz$ be an element of $\st$.  Recall that $a_n = g'(x) \equiv
1 \pmod{p}$ and $b_n = (g(x)-x)/p^n \equiv 0 \pmod{p}$.  Then
\begin{align*}
p^{n+1}b_{n+1} =g(y)-y&=g(x+p^nz)-(x+p^nz) \\
&\equiv g(x)-x + p^n z (g'(x)-1) + p^{2n}z^2 g''(x)/2 \pmod{p^{3n}} \\
&\equiv p^n b_n + p^n z (a_n -1) + p^{2n}z^2 g''(x)/2 \pmod{p^{3n}}
\end{align*}
so
\begin{equation*}
pb_{n+1} \equiv b_n + z(a_n -1) + p^n z^2 g''(x)/2 \pmod{p^{2n}}.
\end{equation*}
Similarly, 
\begin{equation*}
a_{n+1}=g'(y)=g'(x+p^nz) \equiv a_n + p^n zg''(x) \pmod{p^{2n}}.
\end{equation*}
Combining these two expressions gives
\begin{align*}
\frac{z}{2} (a_{n+1}+a_n -2) &\equiv \frac{z}{2} (2a_n + p^n zg''(x) -2)
\pmod{p^{2n}} \\
&\equiv pb_{n+1}-b_n \pmod{p^{2n}}.
\end{align*}

Now we apply this result.  Let $A=\ord(a_n -1)$ and $B=\ord(b_n)$.  We know
that $A,B\geq 1$.  Since $a_{n+1} \equiv a_n \pmod{p^n}$, we have
$\ord(a_{n+1}-1)=A$ if $A<n$, and $\ord(a_{n+1}-1)\geq n$ if $A\geq n$.
Now,
\begin{enumerate}  
\item
If $B<A$ and $B<n$, then $\ord(b_{n+1})=B-1$.
\item
If $A\leq B$ and $A<n$, then $b_n + z(a_n-1) \equiv pb_{n+1}\pmod{p^n}$.
There is a unique $z\pmod{p}$ for which $b_n+z(a_n-1)\equiv 0
\pmod{p^{A+1}}$, so that $\ord(b_{n+1})\ge A$.  For all other $z
\pmod{p}$, $b_n + z(a_n -1)$ is divisible by $p^A$ but not by $p^{A+1}$,
so that $\ord(b_{n+1})=A-1$.
\item
If $A,B\geq n$ then $p^n$ divides $(a_{n+1}-1)$ and $p^{n-1}$ divides
$b_{n+1}$.
\end{enumerate}
Interpreting these results in terms of the tree structure, we see that:
\begin{itemize}
\item
If $B<A$ and $B<n$, then every lift of $\sigma$ splits $(B-1)$ times,
and then grows.
\item
If $A\leq B$ and $A<n$, then every lift of $\sigma$, except for one,
splits $(A-1)$ times, and then grows.  The single exceptional lift
behaves precisely the same way as does $\sigma$.
\item If $A,B\geq n$ then every lift of $\sigma$ splits $n-1$ times, but
we do not know what happens to their descendants at level $2n$.
\end{itemize}

Note that, in addition to using the above results by computing $A$ and
$B$ to predict the structure of the tree, we can use the results by
observing the tree to determine which case we are in.

In cases 1 and 2, the subtrees rooted at $\sigma$ have the following
structures:
\pagebreak
\begin{figure}[h]
\centerline{\psfig{figure=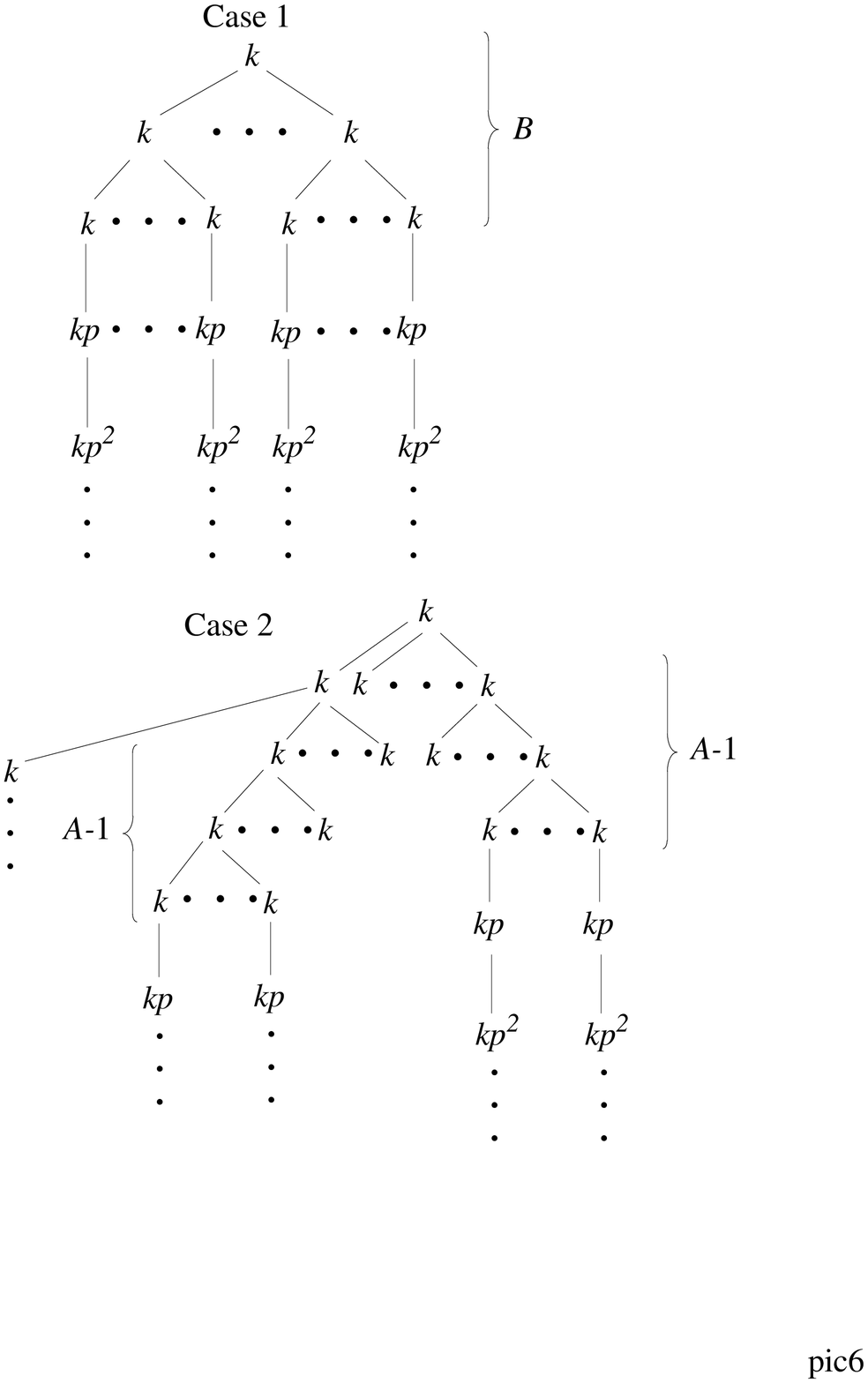,height=6in,width=5in,bblly=2in,bbury=11in,bbllx=0in,bburx=7in,clip=}}
\end{figure}


\subsection{Tails}
\label{tails}

A {\em tail} of $f_n$ is a sequence of elements of $\zpn$ of the form
$y_1,y_2,\dots,y_\ell$, where $y_{j+1}=f_n(y_j)$, and none of the $y_j$
is in the image of $f_n^m$, for $m$ sufficiently large.  All of the
points of $\zpn$ either lie on tails of $f_n$ or in cycles of $f_n$.

Suppose $(x_1,x_2,\dots,x_k)$ is a cycle of $f_1$.  If
$f'(x_i)\not\equiv 0\pmod{p}$, then $f_n$ maps $\{x\mid x\equiv
x_i\pmod{p}\}$ to $\{x\mid x\equiv x_{i+1} \pmod{p}\}$ bijectively.
This follows by induction on $n$.  Let $x\in\zpnmone$ be congruent to
$x_i\pmod{p}$.  Then the $p$ elements of $\zpn$ that are congruent to
$x\pmod{p^{n-1}}$ map bijectively to the $p$ elements of $\zpn$ that are
congruent to $f(x)\pmod{p^{n-1}}$, since 
\begin{equation*}
f(x+p^{n-1}y) \equiv f(x) +
p^{n-1} y f'(x) \equiv f(x) + p^{n-1} y f'(x_i)\pmod{p^n}.
\end{equation*}

Thus, if $f'(x_i)\not\equiv 0\pmod{p}$ for $i=1,\dots,k$, then all of
the elements of $\zpn$ which are congruent to $x_1,\dots,x_k\pmod{p}$
lie on cycles of $f_n$.

However, if $f'(x_i)\equiv 0\pmod{p}$, then $f_n$ maps the $p$ elements
of $\zpn$ that are congruent to $x\pmod{p^{n-1}}$ all to the same
element of $\zpn$, by the above computation.  Thus, the elements of
$\zpn$ which are congruent to $x_1,\dots,x_k\pmod{p}$ contain only a
single cycle of length $k$, and the remaining points lie on tails of
$f_n$.  If $y_1,\dots,y_\ell$ is such a tail, with $y_j\equiv
x_i\pmod{p}$, then $y_{j+1}\pmod{p^2}$ must be on the cycle of $f_2$.
Similarly $y_{j+k+1}\pmod{p^3}$ must be on the cycle of $f_3$, and so
on.  By induction, $y_{j+(n-2)k+1}$ must be on the cycle of $f_n$.
Thus, for such a cycle of $f_n$, the maximum length of a tail leading to
that cycle is $p+(n-2)k$.

Of course, the tails of $f_n$ form trees, with every tail eventually
leading to a cycle, but possibly first joining another tail.  The above
result gives a bound on how long it takes for all the tails to coalesce
into the cycle.

If $f'(x_i)\equiv 0 \pmod{p}$, but $f''(x_i)\not\equiv 0 \pmod{p}$, then
we can describe precisely how $f_n$ maps $\{x\mid x\equiv x_i\pmod{p}\}$
into $\{x\mid x\equiv x_{i+1}\pmod{p}\}$.  The preimages of the points
in the image of $f_n$ have sizes $p^j$ or $2p^j$; precisely, for $1\leq
j < n/2$ there are $p^{n-2j-1}(p-1)/2$ preimages of size $2p^j$, and
there is a single preimage of size $p^{\lfloor n/2\rfloor}$.  The proof
is similar to many we have already presented.


\section{Periodic orbits of $f$}
\label{orbits}

In this section we describe the possible lengths of periodic orbits in
the $p$-adic integers $\zp$ for a polynomial $f(x)\in
\z[x]$.\footnote{Or, more generally, a polynomial in $\zp[x]$.}  Such an
orbit corresponds to a sequence of cycles of $f_n$, for
$n=1,\dots,\infty$, where each cycle is a lift of its predecessor.  The
lengths of the cycles are bounded, and length of the orbit is the lim
sup of the lengths of the cycles.  All the relevant properties of $\zp$
are presented in Appendix~\ref{p-adics}.  We use the term `periodic
orbit' for $\zp$, while we reserve the term `cycle' for $\zpn$.

Let $x\in\zp$ lie in a periodic orbit of $f$ of length $c$.  For each
$n$, let $c_n$ be the length of the cycle $\sigma_n$ of $f_n$ containing
$x\pmod{p^n}$.  Thus, $c_1 \leq c_2 \leq \dots = c$.  Clearly $c_1\leq
p$.  If $\sigma_1$ grows tails, then each $c_n=c_1$, so $c=c_1$.
Otherwise, whenever some $\sigma_n$ either splits or grows, all further
$\sigma_{n+i}$ either split or grow, so $c_{n+i+1}/c_{n+i}$ is always
either 1 or $p$, thus $c/c_n$ is a power of $p$.  If $\sigma_1$
partially splits, then its lifts of length $dc_1$ either split or grow,
and its lift of length $c_1$ partially splits just as does $\sigma_1$.

Thus, there are three possibilities:
\begin{enumerate}
\item $c=c_1$, if either $\sigma_1$ grows tails or every $\sigma_{i+1}
$ is the single $c_i$-cycle lift of $\sigma_i$, which partially splits.
\item $c/c_1$ is a power of $p$, if $\sigma_1$ splits or grows.
\item $c/dc_1$ is a power of $p$, if $\sigma_1$ partially splits but some
$\sigma_{i+1}$ is a $dc_i$-cycle lift of $\sigma_i$ (which partially
splits).
\end{enumerate}

We showed in Section~\ref{grow} that, for $n \geq 2$, whenever $\sigma$
grows, its lift also grows; then that cycle's lift grows, and so on.
So, under the hypothesis that $c=\limsup_{i\to\infty} c_i$ is finite, we
can never have $c_{n+1}=pc_n$ for $n\geq 2$.  Thus, $c=c_1$ or $c=dc_1$
or $c=c_2=pc_1$.

For $p>3$, we show in Appendix~\ref{growth-appendix} that the case
$c=c_2=pc_1$ never occurs.

In summary, any periodic orbit of $f$ in $\zp$ has length at most $p^2$,
and this length is the product of a positive integer not exceeding $p$
and a divisor of $p-1$ (except if $p=3$, in which case length 9 can
occur).  Note that this upper bound on the lengths of periodic orbits
implies, if $f$ is not linear, an upper bound on the number of periodic
points.  For, any element of an orbit of length $c$ must be a root of
the polynomial $f^c(z)-z$, which has only finitely many roots.


\section{Lifts of a periodic orbit}
\label{lifts}

In this section we describe the behavior of cycles which separate from a
periodic orbit of $f$ at some stage.  Precisely, let $\alpha \in \zp$ be
an element of a periodic orbit of $f$ of length $k$, so that
$g(\alpha)=\alpha$ for $g=f^k$, and assume that $g'(\alpha)\not\equiv 0
\pmod{p}$.  Let $c_n=k$ (i.e., $\alpha\pmod{p^n}$ is in a $k$-cycle of
$f_n$, not a shorter cycle), and let $y\in\zp$ have $n=\ord(y-\alpha)$.
Then $\alpha\pmod{p^{n+1}}$ and $y\pmod{p^{n+1}}$ lie in different
cycles of $f_{n+1}$.  We will say that the cycle containing $y$
separates from $\alpha$ at level $n+1$.

Let $d$ be the order of $g'(\alpha)\pmod{p}$.  Then we know that $y
\pmod{p^{n+1}}$ is in a $kd$-cycle of $f_{n+1}$.  Let $h=g^d=f^{kd}$.
Then $h(\alpha)=\alpha$ and $h'(\alpha)\equiv 1\pmod{p}$.

Suppose that $h'(\alpha)=g'(\alpha)^d \neq 1$.  Let
$m=\ord(h'(\alpha)-1)$.  Then
\begin{align*}
h(y)-y &= h(\alpha+(y-\alpha)) -\alpha -(y-\alpha) \\
&\equiv h(\alpha) + (y-\alpha)h'(\alpha) -\alpha -(y-\alpha) \pmod{p^{2n}} \\
&\equiv (y-\alpha) (h'(\alpha)-1) \pmod{p^{2n}}.
\end{align*}
Thus, if $n>m$, $\ord(h(y)-y)=n+m$.\footnote{In the next section, we
will show that the same conclusion holds if $n>m/d$.} In this case,
$y\pmod{p^{n+m}}$ is a fixed point of $h_{n+m}$, but $y\pmod{p^{n+m+1}}$
is not a fixed point of $h_{n+m+1}$.  So $y\pmod{p^{n+m}}$ lies in a
$kd$-cycle of $f_{n+m}$, but $y\pmod{p^{n+m+1}}$ lies in a longer cycle
of $f_{n+m+1}$.  It follows that the $kd$-cycle of $f_{n+m}$ which
contains $y\pmod{p^{n+m}}$, and all of its descendants, must grow.

Note that in the above case there will always be some $n$ such that the
behavior of the infinite subtree consisting of cycles which separate
from $\alpha$ at levels greater than $n$ is determined by the finite
tree up to level $n$.  For there will be some cycle which separates from
$\alpha$ at level $n+1$, and then splits $m$ times where $m<n$.  By the
above argument, all cycles which separate from $\alpha$ at higher levels
will behave the same way.

In fact, by observing only a finite part of the tree, we can determine
that we are in that case.  If we have a cycle which partially splits,
then we know that it has a lift which partially splits, and so on, so
each point on that cycle corresponds to a periodic element $\alpha$.
Then, if one of the other cycles which is a lift of that cycle behaves
as above (splits $m<n$ times, then grows), then we know that we are in
the above case and all cycles which separate from $\alpha$ at higher
levels will behave the same way.

If we have a periodic element $\alpha$ which is on a cycle which splits
completely, and it has a lift which splits $m<n$ times and then grows,
then simply by observing that feature of the tree, by the results of
Section~\ref{split} we must be in case~2 of that section, and so we know
that there is a periodic point $\alpha$ on the cycle, and the above
results apply.

However, when $h'(\alpha)=1$, it does not seem to be true that by
observing a finite part of the tree we can predict all subsequent
behavior, nor can we determine that $h'(\alpha)=1$ by observing only a
finite portion of the tree.

Suppose that $h'(\alpha)=1$.  If $f$ is not linear, then $h'$ is not
constant, so there is an integer $\ell\geq 2$ for which
$h^{(\ell)}(\alpha) \neq 0$ while $0=h^{(2)}(\alpha) =\dots =h^{ (\ell
-1)}(\alpha).$\footnote{In the next section, we will show that
$\ell>d$.} Let $m=\ord(h^{(\ell)}(\alpha)/\ell !)$.  Then
\begin{align*}
h(y)-y &=h(\alpha+(y-\alpha)) -\alpha - (y-\alpha) \\
&=h(\alpha) + (y-\alpha)h'(\alpha) +\dots -\alpha-(y-\alpha) \\
&= (y-\alpha)^\ell h^{(\ell)}(\alpha)/\ell ! + (y-\alpha)^{\ell +1}
h^{(\ell +1)}(\alpha)/(\ell +1)! + \dots \\
&\equiv (y-\alpha)^\ell h^{(\ell)}(\alpha)/\ell ! \pmod{p^{n(\ell +1)}}.
\end{align*}
Thus, if $n>m$, $\ord(h(y)-y)=n\ell+m$.  In this case, the image of
$y\pmod{p^{n\ell+m}}$ is in a $kd$-cycle of $f_{n\ell+m}$ which grows
and all of whose descendants grow.  Thus, the lifts of $\alpha$ which
separate from it at any stage $n+1$, where $n>m$, will split
$n(\ell-1)+(m-1)$ times and then grow.  (And, since this is greater than
or equal to $n$, the above results for the case $h'(\alpha)\neq 1$ never
apply.) The lifts of $\alpha$ which separate from it at stage $n+1$,
where $n\leq m$, will split at least $n\ell-1$ times, but we do not know
whether they then grow.

We consider the case $h'(\alpha)=1$ to be pathological; it did not arise
in any of thousands of random examples we studied.  We can construct an
example, though: let $p=3$ and $f(x)=x+3x^2$.  Then $f(0)=0$, so take
$\alpha=0$.  Since $f(3^n\beta)\equiv 3^n\beta \pmod{3^{2n+1}}$ for any
$n,\beta$, the cycles which separate from 0 at level $n+1$ split $n$
times and then grow.


\section{Improving the bounds}
\label{improvements}

In the previous section we described the dynamics of $f$ sufficiently close
to a periodic orbit; in this section we will show that the same results hold
in a somewhat larger neighborhood of the periodic orbit.  
We will prove the following result:

\begin{proposition}
\label{ladeda}
If $f(x)\in \zp[x]$ has a periodic orbit of length $k$ containing
$\alpha \in \zp$, $d>1$ is the order of $(f^k)'(\alpha) \pmod{p}$, and
$h=f^{kd}$, then each of $h^{(2)}(\alpha),\dots,h^{(d)}(\alpha)$ is
divisible by $(h'(\alpha)-1)$, in $\zp$.
\end{proposition}

Our interest is in the following two corollaries:

\begin{corollary}
\label{cor1}
Under the hypotheses of Proposition~\ref{ladeda}, if $h'(\alpha)=1$ then 
$h^{(2)}(\alpha)=\dots=h^{(d)}(\alpha)=0$.
\end{corollary}

\begin{corollary}
\label{cor2}
Under the hypotheses of Proposition~\ref{ladeda}, if $m={\rm ord}_p(h'(\alpha)-1)$
and $y\in\zp$ has $n={\rm ord}_p(y-\alpha)$, then
\begin{equation*}
h(y)-y \equiv (y-\alpha)(h'(\alpha)-1)
\pmod{p^{\min\{n(d+1),2n+m\}}}.
\end{equation*}
\end{corollary}

\begin{proof}
Observe that 
\begin{equation*}
h(y)-y = h(\alpha)-\alpha + (y-\alpha)(h'(\alpha)-1) +
(y-\alpha)^2h''(\alpha)/2!+\cdots,
\end{equation*}
where $h(\alpha)=\alpha$ and $p^{ni+m}$ divides $(y-\alpha)^ih^{(i)}
(\alpha)/i!$ for $2\leq i\leq d$.
\end{proof}

The following corollary was used in Section~\ref{partial}:

\begin{corollary}
\label{useful}
Under the hypotheses of Corollary~\ref{cor2}, if we define $a_{n+1} =
h'(y)$, then 
\begin{equation*}
\min\{{\rm ord}_p(h(y)-y)-n,nd\} = \min\{{\rm ord}_p(a_{n+1}-1),nd\}.
\end{equation*}
\end{corollary}

\begin{proof}
We have 
\begin{equation*}
a_{n+1}=h'(y)=h'(\alpha)+(y-\alpha)h''(\alpha)+\cdots.
\end{equation*}
But $(y-\alpha)^{i-1}h^{(i)}(\alpha)$ is divisible by $p^{m+n(i-1)}$ for
$2\leq i \leq d$, so 
\begin{equation*}
a_{n+1} \equiv h'(\alpha) \pmod{p^{\min\{n+m,nd\}}}.
\end{equation*}
Thus, $\min\{\ord(a_{n+1}-1),n+m,nd\}=\min\{m,n+m,nd\}$.  From
Corollary~\ref{cor2}, 
\begin{align*}
\min\{\ord(h(y)-y)-n,n+m,nd\} &= \min\{m,n+m,nd\} \\
&= \min\{\ord(a_{n+1}-1),n+m,nd\}.
\end{align*}
But $n+m>\ord(a_{n+1}-1)$, so the minimum of the right-hand side is less
than $n+m$, so the minimum of the left-hand side is also less than
$n+m$.
\end{proof}

\begin{proof}[Proof of Proposition~\ref{ladeda}]
First, we will translate the
periodic orbit so that it passes through 0; this simplifies the algebra
in our proof.  Let $T:x\mapsto x+\alpha$.  Then the function $\hat f =
T^{-1}fT$ also has a periodic orbit of length $k$, namely $(0,
f(\alpha)-\alpha,f^2(\alpha)-\alpha,\dots,f^{k-1}(\alpha)-\alpha)$.
But $\hat f(x) = f(x+\alpha)-\alpha$, so $\hat f'(x)=f'(x+\alpha)$, and
similarly for higher derivatives; likewise, any iterate $\hat f^{\ell}=
T^{-1}f^{\ell}T$, so $(\hat f^{\ell})^{(t)}(x)=
(f^{\ell})^{(t)}(x+\alpha)$.  So computations assuming that $\alpha=0$
will also hold for arbitrary $\alpha$.

Now, for $g=f^k$ we have $g(x)=g'(0)x+\OO(x^2)\in\zp[x]$ and
$h(x)=g^d(x)=g'(0)^dx +h_2x^2+h_3x^3+\dots+h_dx^d+\OO(x^{d+1})$, where
$\OO(x^j)$ denotes a polynomial in $x$ in which every term has degree at
least $j$.  Since $d$ is the order of $g'(0)\pmod{p}$, each of
$g'(0)-1,\dots,g'(0)^{d-1}-1$ is coprime to $p$.  Thus, for any $\ell
\leq m=\ord(g'(0)^d-1)$, we can project to $R=\zp/p^{\ell}\zp$ and
apply the following lemma, which implies that each of $h_2,\dots,h_d$ is
divisible by $p^{\ell}$, completing the proof of the Proposition.
\end{proof}

\begin{lemma}
For any commutative ring $R$ and any primitive $d^{\rm th}$ root of unity
$\zeta\in R$ such that none of $\zeta-1, \zeta^2-1,\dots, \zeta^{d-1}-1$
is a zero-divisor, let $g(x)=\zeta x+\OO(x^2)\in R[x]$ and
$h(x)=g^d(x)=x+ax^i+\OO(x^{i+1})$, where $a\neq 0$ is the first nonzero
coefficient of $h(x)$ of degree greater than 1.  Then $i\equiv 1
\pmod{d}$, and in particular $i\geq d+1$.
\end{lemma}
\begin{proof}
Write $g(x)=g_1 x +g_2 x^2 +\dots$.  Then the compositions
\begin{equation*}
h\mycirc g = g_1 x + \dots + g_{i-1}x^{i-1} + (g_i+a\zeta^i) x^i + \OO(
x^{i+1})
\end{equation*}
and 
\begin{equation*}
 g\mycirc h = g_1 x + \dots + g_{i-1}x^{i-1} + (g_i+\zeta a) x^i + 
\OO(x^{i+1}).
\end{equation*}
Since $g\mycirc h = h\mycirc g$, the coefficients of $x^i$ are equal, so
$a\zeta^i = a\zeta$, so $a\zeta(\zeta^{i-1}-1)=0$.  Since $\zeta$ is a root of
unity, it is not a zero-divisor.  Therefore if $i\not\equiv 1
\pmod{d}$, then $\zeta^{i-1}-1$ would be a zero-divisor, contradicting
the hypothesis.
\end{proof}


\section{Odds and ends}


\subsection{Analyzing polynomials}

For a given polynomial, our results generally allow us to find the cycle
structure of $f_n$ rather quickly.  We can compute the first few levels
of the tree directly, and then our results will usually imply the
structure of the entire tree.  We have done this for thousands of
randomly selected polynomials, for small primes $p$; in theory one
should be able to construct polynomials which will take us arbitrarily
long to analyze, but these polynomials seem to be extremely rare.  Also,
the numbers $a_n$ and $b_n$ are sometimes useful for determining the
structure of the remainder of the tree.

The tree shown in Section~\ref{philosophy} is a typical example.  For
this tree, levels 0--3 suffice to determine the structure of the entire
tree.  Once we observe a cycle $(\text{mod } {p^2})$ which splits, for which one
lift grows and another splits, we know that the lift which splits will
behave in the same way.  Also, for the 1-cycle $(\text{mod } {p})$ which
partially splits, since its 2-cycle lift splits 0 times before growing,
and $0<kr-1=1$, this behavior must persist.

Conversely, our results also allow us to construct polynomials with
desired cycle structures $(\text{mod } {p^n})$.  For instance, we can construct
polynomials having periodic orbits (in $\zp$) of length $kr$, for any
$1\leq k \leq p$ and any $r$ dividing $p-1$.


\subsection{Polynomial with 3-adic 9-cycle}

We mention the polynomial $f(x)=2+x+3x^2+x^3+3x^4+2x^5$, which has a
3-adic 9-cycle, since it has a 9-cycle (mod 81) for which $\ord(a_4
-1)=3$ and $\ord(b_4)=4$, namely the cycle containing 0 (mod 81).

A $p$-adic cycle of length $p^2$ is impossible for $p>3$, by the results
of Appendix~\ref{growth-appendix}.


\subsection{Permutation polynomials and single-cycle polynomials}

We give a straightforward method for determining whether a given
polynomial $f(x)\in\z[x]$ induces a permutation
$f_n:\zpn\rightarrow\zpn$.  We claim that, for $n\geq 2$, $f_n$ is a
permutation if and only if $f_1$ is a permutation and $f'(x)$ has no
roots in $\zpp$; it follows that, for any $n\geq 2$, $f_n$ is a
permutation if and only if $f_2$ is a permutation.

We prove the claimed result by induction.  Certainly, if $f_n$ is a
permutation, then $f_{n-1}$ is a permutation, which implies that $f_1$
is a permutation.  Now, given that $f_{n-1}$ is a permutation, $f_n$
will be a permutation if and only if, for each integer $x$, the numbers
$f(x),f(x+p^{n-1}),\dots,f(x+p^{n-1} (p-1))$ are all distinct
$(\text{mod } {p^n})$; but $f(x+p^{n-1}t)\equiv f(x)+p^{n-1}t f'(x) \pmod{p^n}$,
so $f_n$ is a permutation if and only if $f'$ has no roots in $\zpp$.
This completes the proof.

We can also give a simple criterion for when $f_n$ is a single cycle of
length $p^n$.  For $p>3$, for any $n\geq 2$, $f_n$ is a $p^n$-cycle if
and only if $f_2$ is a $p^2$-cycle.  For $p=3$, for any $n\geq 3$, $f_n$
is a $3^n$-cycle if and only if $f_3$ is a $3^3$-cycle.


\section{Further generality}
\label{etc.}

There are more general situations in which our arguments, perhaps with slight
modifications, will apply.  They include various combinations of the following:
\begin{itemize}
\item First of all, we can replace $\z$ by $\zp$ in all of our arguments.
\item Most of the results which we have derived for polynomials also hold for
rational functions whose denominators have no roots in $\zpp$; we will
show this in Appendix~\ref{rational-functions}.
\item More generally, we can consider rational functions over $\qp$ having
``good reduction'' (mod $p$) at all points of a cycle in $\PP^1(\qp)$. 
\item All of the above proofs work just as well for power series
over $\zp$, with one caveat: if the power series only converges on $p\zp$,
then we must only consider elements of $p\zp$, and in particular $f(0)$ must
be divisible by $p$.
\item The above arguments apply, in modified form, if we replace $\zp$
by the valuation ring of any finite extension of $\qp$.
\item For polynomials with coefficients in a number field, we can pick a
good prime of the number field (almost any would do) and apply the
results for the valuation ring of the completion of the field at that
prime, to give bounds on the cycle lengths.
\item Our basic approach yields interesting results for polynomial
mappings from $\z^n$ to $\z^n$ (thanks to Greg Kuperberg for pointing this 
out).
\end{itemize}
We have studied all of the above, and we have numerous partial results;
we hope eventually to write a comprehensive paper covering at least the
above situations.


\newpage
\appendix

\section{A quick introduction to $p$-adics.}
\label{p-adics}

For the reader's convenience, we set forth the basic properties of
$p$-adic integers which we use in Sections~\ref{orbits}, \ref{lifts},
\&~\ref{improvements}.  The $p$-adic integers are the projective limit of
the rings $\zpn$.  Explicitly, an element of $\zp$ is a sequence
$(x_1,x_2,\dots)$, where $x_n \in \zpn$, such that $x_m \equiv x_n
\pmod{p^n}$ for all $m>n$.  Addition and multiplication are defined
component-by-component, which makes $\zp$ into a ring.  Note that $\zp$
contains $\z$, since any nonnegative integer $n$ is represented by
$(n,n,n,\dots)$.  Also note that $\zp$ is a domain, namely, $a\cdot b=0$
only happens when $a$ or $b$ is zero.  It makes sense to reduce elements
of $\zp$ modulo $p^n$, in the usual ring-theoretic way or just by
extracting the $n\tth$ component; the ring $\zp/{p^n}\zp$ is canonically
isomorphic to $\zpn$.  Finally, the invertible elements of $\zp$ are
precisely the elements not divisible by $p$ ({\em i.e.}, the elements
for which $x_1 \not\equiv 0 \pmod{p}$), for one can write down the
inverse of such an element component-by-component, in much the same way
as one multiplies in $\zp$.

The $p$-adic integers can be visualized as an infinite tree much like
the ones we have described above.  Construct the infinite $p$-ary tree,
analogous to the infinite binary tree for $p=2$, viewing the nodes on
the $n\tth$ level as classes $(\text{mod } {p^n})$, labeled in such a way that a
node $(\text{mod } {p^n})$ is connected by an edge to the $p$ nodes
$(\text{mod } p^{n+1})$ which are congruent to the first class $(\text{mod } {p^n})$.
The first few levels of the tree for $p=3$ are pictured below:

\begin{figure}[h]
\centerline{\psfig{figure=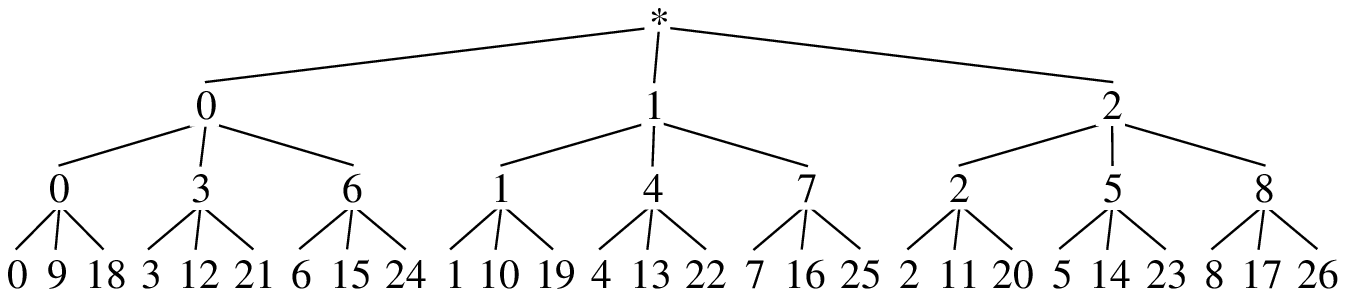,height=1.6in,width=5in,bblly=1in,bbury=3in,bbllx=1in,bburx=7in,clip=}}
\end{figure}

Then the $p$-adic integer $(x_1,x_2,\dots)$ corresponds to the infinite
path down the tree which passes through each node $x_n \pmod{p^n}$.
$\zp$ is the set of all such infinite paths in the tree.  This
interpretation makes it easy to see when two $p$-adic integers are
congruent $(\text{mod } {p^n})$: if and only if their paths coincide for the
first $n$ levels.  If $\ord(x-y)=n$, then the paths $x$ and $y$ coincide
at the first $n$ levels, and separate at level $n+1$.


\newpage
\section{Cycle growth}
\label{growth-appendix}

We prove that, for $p>3$, if a cycle of $f_1$ grows, then its lift also
grows.  From previous results we then know that the next lift grows, and
so on.  We also describe when this fails for $p=3$.

Say our polynomial $f$ has a $k$-cycle mod $p$ which contains $x$.  Let
$g=f^k$, which has $x$ as a fixed point mod $p$, but for which $x$ is in
a $p$-cycle mod $p^2$.  We must show that $x$ is in a $p^2k$-cycle of
$f_3$, or equivalently that $x$ is in a $p^2$-cycle of $g_3$.  Let
$a=a_1=g'(x)$, $b=b_1=(g(x)-x)/p$, and $c=g''(x)/2$.  Then $a\equiv 1
\pmod{p}$ and $b \not\equiv 0 \pmod{p}$.

First, we show that, for each $i\geq 1$,
\begin{equation*}
g^i(x) \equiv x+pb\sum_{j=0}^{i-1}a^j + p^2cb^2
\sum_{j=0}^{i-2}a^{i-2-j} (1+a+\dots +a^j)^2 \pmod{p^3}.
\end{equation*}
For, this is
true for $i=1$, and inductively
\begin{align*}
g^{i+1}(x) &\equiv 
g\left(x+p\left(b\sum_{j=0}^{i-1}a^j + pcb^2 \sum_{j=0}^{i-2}
a^{i-2-j}(1+a+\dots +a^j)^2\right)\right) \pmod{p^3} \\
&\equiv g(x) +pb\sum_{j=1}^i a^j +p^2cb^2 \sum_{j=0}^{i-2}a^{(i+1)-2
-j}(1+a+\dots +a^j)^2 \\
&+\, p^2 cb^2 \left(\sum_{j=0}^{i-1}a^j\right)^2 \pmod{p^3} \\
&=x +pb\sum_{j=0}^i a^j +p^2 cb^2 \sum_{j=0}^{i-1} a^{(i+1)-2-j}(1+a+\dots+
a^j)^2,
\end{align*}
which completes the induction.

Now,
\begin{align*}
g^p(x) &\equiv x+pb \sum_{j=0}^{p-1}a^j +p^2cb^2
\sum_{j=0}^{p-2}a^{p-2-j} \left(\sum_{\ell =0}^j a^\ell\right)^2 \pmod{p^3} \\
&\equiv x+pb\sum_{j=0}^{p-1}a^j +p^2 cb^2 \sum_{j=0}^{p-2}(j+1)^2
\pmod{p^3}
\,\,\,\text{ (since $a\equiv 1$ (mod $p$))} \\
&=x+pb\sum_{j=0}^{p-1}a^j +p^2 cb^2 \frac{(p-1)(p)(2p-1)}{6}
\end{align*}
and, for $p>3$, the last term is $0 \pmod{p^3}$, so 
\begin{equation*}
g^p (x) \equiv x+pb\sum_{j=0}^{p-1}a^j \pmod{p^3}.
\end{equation*}
But, as shown in Section~\ref{grow},
$\sum_{j=0}^{p-1}a^j$ is not divisible by $p^2$.  Thus, $g^p(x) 
\not\equiv x \pmod{p^3}$, so the $p$-cycle of $g_2$ which includes $x$
does not split, hence it grows.

Using the above methods, we can describe when a cycle $(\text{mod } 3)$ will
grow and then split.  For $p=3$ we have
\begin{align*}
g^p(x) &\equiv  x+pb\sum_{j=0}^{p-1}a^j +p^2 cb^2 \frac{(p-1)(p)(2p-1)}{6} 
\pmod{p^3} \\
&\equiv  x+pb\sum_{j=0}^{p-1}a^j - p^2 c \pmod{p^3} \\
&\equiv  x+p^2 b - p^2 c \pmod{p^3},
\end{align*}
so a cycle (mod~3) which grows will then split if and only if $b\equiv
c\pmod{3}$.


\newpage
\section{Rational functions}
\label{rational-functions}

Let $h=f/g$ be a ratio of polynomials $f,g\in\z[x]$ such that $g$ takes values
coprime to $p$ on any cycle being considered; in particular, this condition
certainly holds if $g$ has no roots in the field $\zpp$.  We will show that
the results we have derived for polynomials over $\z$ also hold for $h$.
We do this by constructing a sequence of polynomials $h_n\in\z[x]$ such that
$h$ and $h_n$ agree $(\text{mod } {p^n})$ on the cycles being considered, and the 
$a_i$'s and $b_i$'s for the various $h_n$ are compatible.

Precisely, put 
\begin{equation*}
h_n(x)=f(x)\cdot g(x)^{\phi(p^{2n})-1},
\end{equation*}
where $\phi$ is the Euler quotient function; then $h_n(x)\in\z[x]$.
For any $x$ such that $p\nmid g(x)$, $h_n(x)\equiv h(x)\pmod{p^{2n}}$.
 Let $a_{i,n}$ and
$b_{i,n}$ be the values of $a_i$ and $b_i$ for the polynomial $h_n$, for
$i\leq n$, and 
say that $x$ is in a cycle of $h_n \pmod{p^i}$ of length $\alpha_{i,n}$.
Then 
\begin{equation*}
b_{i,n}=\frac{h_n^{\alpha_{i,n}}(x)-x}{p^i}\equiv 
\frac{h^{\alpha_{i,n}}(x)-x}{p^i} \pmod{p^n}.
\end{equation*}
Next,
\begin{align*}
h_n'(x)&=f'(x)\cdot g(x)^{\phi(p^{2n})-1}+f(x)\cdot \bigl(\phi(p^{2n})-1
\bigr)\cdot g(x)^{\phi(p^{2n})-2}\cdot g'(x) \\
&\equiv\frac{f'(x)}{g(x)}+f(x)\cdot(p^{2n}-p^{2n-1}-1)\cdot\frac{g'(x)}{g(x)^2
} \pmod{p^{2n}}\\
&\equiv \frac{g(x)f'(x)-f(x)g'(x)}{g(x)^2} \pmod{p^{2n-1}} \\
&=h'(x),
\end{align*}
so
\begin{align*}
a_{i,n} &=(h_n^{\alpha_{i,n}})'(x) = \prod_{\ell=0}^{\alpha_{i,n}-1}h_n'(
h_n^\ell(x)) \equiv \prod_{\ell=0}^{\alpha_{i,n}-1}h'(h^\ell(x)) \pmod{p^{2n-1}}
\\
&=(h^{\alpha_{i,n}})'(x).
\end{align*}

Now, for any $n>i$, $h$ agrees with $h_n \pmod{p^i}$; thus, each $\alpha_{i,n}$
equals the length of the cycle of $h$ mod $p^n$ containing $x$.  Hence,
the classes $a_{i,n} \pmod{p^i}$ and $b_{i,n}\pmod{p^i}$ are independent of $n$.
This shows the compatibility of the $h_n$; thus, because our earlier
results apply to each $h_n$, they apply as well to the function $h$.


\end{document}